\documentclass[12pt]{article}

\usepackage{url}
\usepackage{pslatex}
\usepackage{amsmath,amsthm,amssymb,amsfonts, stmaryrd}
\usepackage{fancybox}
\usepackage{color, graphicx}
\usepackage{hyperref}
\usepackage[T1]{fontenc} 
\usepackage[usenames,dvipsnames,svgnames,table]{xcolor}
\usepackage{framed}
\newtheorem{theorem}{Theorem}[section]

\newtheorem{corollary}[theorem]{Corollary}

\newtheorem{definition}[theorem]{Definition}
\newtheorem{example}[theorem]{Example}

\newtheorem{lemma}[theorem]{Lemma}
\newtheorem{notation}[theorem]{Notation}

\newtheorem{proposition}[theorem]{Proposition}

\newcommand{\vs}[1]{\langle #1 \rangle}

\newcommand{\F}{\mathcal{F}}
\newcommand{\G}{\mathcal{G}}
\newcommand{\M}{\mathcal{M}}
\newcommand{\N}{\mathbb{N}}
\newcommand{\Z}{\mathbb{Z}}

\newcommand{\hatpi}[1]{\hat{\pi}{#1}}

\DeclareMathOperator{\inm}{in}

\DeclareMathOperator{\lc}{lc}
\DeclareMathOperator{\lt}{lt}

\begin{document}

\title{Optimal Bounds on the Growth of Iterated Sumsets in Abelian Semigroups}

\author{Shalom Eliahou\footnote{LMPA-ULCO, Calais, France. Email: eliahou(at)univ-littoral.fr}~ and 
Eshita Mazumdar\footnote{Ahmedabad University, India. Email: eshita.mazumdar(at)ahduni.edu.in}}

\date{}

\maketitle

\begin{abstract} We provide optimal upper bounds on the growth of iterated sumsets $hA=A+\dots+A$ for finite subsets $A$ of abelian semigroups. More precisely, we show that the new upper bounds recently derived from Macaulay's theorem in commutative algebra are best possible, i.e., are actually reached by suitable subsets of suitable abelian semigroups. Our constructions, in a multiplicative setting, are based on certain specific monomial ideals in polynomial algebras and on their deformation into appropriate binomial ideals via Gr\"{o}bner bases.
\end{abstract}


\section{Introduction}

Let $A$ be a nonempty finite subset of an abelian semigroup $(G,+)$. Estimating the growth of the iterated sumsets $hA=\underbrace{A+\cdots+A}_h$ as $h$ increases is a core problem in additive combinatorics. 
Khovanskii~\cite{K1,K2} showed that $|hA|$ is asymptotically polynomial in $h$. See also~\cite{N2, NR}. But not much is known about this polynomial and, for $h$ small, the behavior of $|hA|$ may wildly vary with $A$, even when $|A|$ is fixed. 
A classical estimate, originally derived using graph theory, is given by Pl\"unnecke's inequality~\cite{Pl}, namely
\begin{equation}\label{plunnecke}
|hA| \le |iA|^{h/i}
\end{equation}
for all $1 \le i \le h$. See~\cite{GR,N1,Pe, TV} for in-depth treatments of this and related inequalities. We recently improved~\eqref{plunnecke} by deriving it from Macaulay's 1927 theorem on the growth of Hilbert functions of standard graded algebras~\cite{EM}. Macaulay's theorem involves a certain operation $a \mapsto a^{\vs{h}}$ on positive integers related to \emph{binomial representations}. In short, if 
$\displaystyle a = \sum_{j=1}^h \binom{a_j}j$ with decreasing integers $a_h > \cdots > a_1 \ge 0$, then $\displaystyle a^{\vs{h}} = \sum_{j=1}^h \binom{a_j+1}{j+1}$, and this is well-defined. See Section 2 for more details. Using this notation, here is part of our improvement to~\eqref{plunnecke} obtained in~\cite{EM}.

\begin{theorem}\label{thm EM} Let $A$ be a nonempty finite subset of an abelian semigroup $G$. Set $d_h=|hA|$ for all $h$. Then $d_0=1$ and
\begin{equation}\label{improvement}
d_{h+1} \le d_h^{\vs{h}}
\end{equation}
for all $h \ge 1$.
\end{theorem}

\begin{example}
For comparison purposes, let $A \subset \Z$ be a subset such that $|6A|=1000$. While Pl\"unnecke's inequality~\eqref{plunnecke} yields 
$$
|5A| \ge 317, \quad |7A| \le 3162,
$$
inequality~\eqref{improvement} yields the much sharper -- and nearly optimal --  bounds
\begin{equation}\label{example sharper}
|5A| \ge 511, \quad |7A| \le 1827.
\end{equation}
See Example~\ref{ex 1000 6} below for the derivation of $|7A| \le 1827$  from $|6A| = 1000$ and~\eqref{improvement}.
\end{example}

Our purpose in this paper is to prove that the upper bounds in Theorem~\ref{thm EM} are best possible. That is, if $(d_i)_{i \ge 0}$ is any sequence of positive integers such that $d_0=1$ and $d_{i+1} \le d_i^{\vs{i}}$ for all $i \ge 1$, then there exists an abelian semigroup $G$ and a subset $A \subseteq G$ such that
\begin{equation}\label{equality}
d_h = |hA|
\end{equation}
for all $h \ge 0$. Our construction of such a pair $G,A$ is in multiplicative notation and proceeds as follows. Let $n=d_1$ and $S=K[X_1,\dots,X_n]$, the $n$-variable polynomial algebra over a field $K$ with its standard grading. Then $G$ will be a multiplicative submonoid of a quotient ring $R = S/J$, where $J$ is an appropriate graded ideal of $S$. Denoting by $\pi \colon S \to R$ the quotient map, and setting $x_j=\pi(X_j)$ for $1 \le j \le n$, we consider the subset $A=\{x_1,\dots,x_n\}$ of $R$ and its $h$-fold iterated product sets $A^h=A\cdots A$. The problem then amounts to uncover a suitable ideal $J$ of $S$ so as to realize, for this subset $A$ of $S/J$, the equality $d_h=|A^h|$ for all $h$. For an almost sharp realization, a specific monomial ideal $J=L$ establishing the converse part of Macaulay's theorem suffices. A sharp realization is then achieved by deforming $L$ into a binomial ideal $\hat{L}$ via a Gr\"{o}bner basis construction so as to preserve the Hilbert function of $S/L$.

\smallskip
The contents of this paper are as follows. Section~\ref{sec background} provides some background on binomial representations, Macaulay's theorem and lexideals. In Section~\ref{sec almost}, after recalling basic facts about monomial ideals, we prove that the bounds provided by Theorem~\ref{thm EM} are almost sharp in an appropriate sense. In~Section~\ref{sec full}, after recalling basic facts about Gr\"{o}bner bases, we proceed to prove the full sharpness of these bounds. The analogous problem restricted to abelian groups remains open. This is briefly discussed in the concluding Section~\ref{sec conclusion}.

\section{Background}\label{sec background}

Given sets $A,B$ in an abelian semigroup $(G,+)$, their \emph{sumset} is $A+B=\{a+b \mid a \in A, b \in B\}$. For $A=B$, we denote $2A=A+A$, and more generally $hA = \underbrace{A+\dots+A}_{h}$ for all $h \ge 2$. Macaulay's theorem involves a certain operation $a \mapsto a^{\vs{h}}$ on $\N$ related to binomial representations, which we now recall.

\subsection{Binomial representation}

\begin{proposition}\label{prop binrep} Let $h \ge 1$ be a fixed integer. Then for any integer $a \ge 1$, there are unique integers $a_h > a_{h-1} > \cdots > a_1 \ge 0$ such that 
$$
a = \sum_{j=1}^h \binom{a_j}j. 
$$
\end{proposition}
\begin{proof} See e.g. the relevant chapters in \cite{BH, HH, P}.
\end{proof}

This expression is called the $h$-\emph{binomial representation of $a$}. Producing it is computationally straightforward: take for $a_h$ the largest integer such that $\binom{a_h}{h} \le a$, and complete that first summand by adding to it the $(h-1)$-binomial representation of $a-\binom{a_h}{h}$. The unicity follows from the classical formula
\begin{equation}\label{classical}
\binom{n+h}h = \sum_{j=0}^{h} \binom{n-1+j}j. 
\end{equation}

\begin{notation}
 Let $a \ge h \ge 1$ be integers. Let
$ \displaystyle
a = \sum_{j=1}^h \binom{a_j}j
$
be its unique $h$-binomial representation. We then denote
$\displaystyle a^{\vs{h}} = \sum_{j=1}^h \binom{a_j+1}{j+1}$. We also set $0^{\vs{h}}=0$.
\end{notation}

Note that the right-hand side $\sum_{j=1}^h \binom{a_j+1}{j+1}$ is a valid $(h+1)$-binomial representation of some positive integer, namely of the integer it sums to.

\begin{example}\label{ex 1000 6} Let $h=6$ and $a=1000$. Then
\begin{align*}
1000 &= \binom{12}{6} +\binom{8}{5} +\binom{6}{4} +\binom{4}{3} +\binom{2}{2} +\binom{0}{1}, \text{whence}\\ 
1000^{\vs{6}} &= \binom{13}{7} +\binom{9}{6} +\binom{7}{5} +\binom{5}{4} +\binom{3}{3}  +\binom{1}{2} = 1827.
\end{align*} This explains the upper bound in~\eqref{example sharper} using 
Theorem~\ref{thm EM}.
\end{example}

\subsection{Macaulay's theorem}
\medskip
Let $R=\oplus_{i \ge 0}R_i$ be a standard graded algebra over a field $R_0=K$. That is, $R$ is a graded commutative algebra which is finitely generated by $R_1$ as a $K$-algebra. It follows that $R_i=R_1^i$, the $i$-fold product set of $R_1$, and that $R_i$ is finite-dimensional as a vector space over $K$ for all $i \ge 0$. The \emph{Hilbert function} of $R$ is the numerical function $i \mapsto d_i=\dim_{K}R_i$.

Macaulay's classical theorem gives necessary and sufficient conditions for any numerical function $i \mapsto d_i$ to be the Hilbert function of a standard graded algebra~\cite{M}. Here it is.

\begin{theorem}[Macaulay]\label{thm macaulay} Let $R = \oplus_{i \ge 0} R_i$ be a standard graded algebra over a field $K$, with Hilbert function $d_i = \dim R_i$. Then $d_0=1$ and
\begin{equation}\label{macaulay}
d_{i+1} \ \le \ d_i^{\vs{i}}
\end{equation}
 for all $i \ge 1$. Conversely, let $(d_i)_{i \ge 0}$ be a sequence of nonnegative integers such that $d_0=1$ and $d_{i+1} \ \le \ d_i^{\vs{i}}$ for all $i \ge 1$. Then there exists a standard graded $K$-algebra $R=\oplus_{i \ge 0} R_i$ such that $d_i=\dim R_i$ for all $i \ge 0$. 
\end{theorem}

With the notation of Theorem~\ref{thm macaulay}, note that if $d_i = 0$ for some $i \ge 2$, then $d_j = 0$ for all $j \ge i$, and this occurs if and only if $R$ is finite-dimensional as a $K$-vector space. A more detailed version of the converse statement, needed for our present purposes, is given below.

\subsection{Lexideals}\label{sec lexideals}

For the converse part in Theorem~\ref{thm macaulay}, the desired algebra $R$ may be constructed as a quotient of a polynomial algebra by a suitable monomial ideal (see Section~\ref{sec monomial ideals}), and more specifically by a \emph{lexideal} $L$. Here are some details needed in the sequel.

Let $(d_i)_{i \ge 0}$ be a sequence of nonnegative integers such that $d_0=1$ and $d_{i+1} \ \le \ d_i^{\vs{i}}$ for all $i \ge 1$. Set $d_1=n$. In the polynomial algebra $S=K[X_1,\dots,X_n]$ over the field $K$, with its standard grading given by $\deg(X_j)=1$ for all $j$, we endow the set $\M$ of monomials in $S$ with the \emph{graded lexicographic order} relative to $X_1 > \dots > X_n$. That is, for $u=\prod_j X_j^{a_j}, v=\prod_j X_j^{b_j} \in \M$, we set $u > v$ if either $\deg(u)>\deg(v)$, or else $\deg(u)=\deg(v)$ and $u$ comes before $v$ lexicographically, i.e. the first nonzero difference $a_j-b_j$ is positive.

\begin{example} With this ordering, the monomials of degree 2 in $K[X_1,X_2,X_3]$ are ordered as follows: 
$$
X_1^2 > X_1X_2 > X_1X_3 > X_2^2 > X_2X_3 > X_3^2.
$$
\end{example}

For all $i \ge 0$, we denote by $\M_i$ the set of monomials of degree $i$ in $S$. Thus $\M_0=\{1\}$, $\M_1=\{X_1,\dots,X_n\}$ and $\M_i=\M_1^i$, the $i$-fold product set of $\M_1$.

\begin{definition} A \emph{lexsegment} is a subset $C$ of $\M_i$ for some $i \ge 1$ such that $C=\{u \in \M_i \mid u \ge v\}$ for some $v \in \M_i$. A \emph{lexideal} $L$ in $S$ is a monomial ideal such that $L \cap \M_i$ is a lexsegment for all $i \ge 1$ such that $L \cap \M_i\not= \emptyset$.
\end{definition}

It is easy to verify that if $C \subseteq \M_i$ is a lexsegment, then $\M_1 C\subseteq \M_{i+1}$ is a lexsegment as well, where $\M_1 C = \{X_j u \mid u \in C, 1 \le j \le n\}$. The converse in Macaulay's theorem may be expressed in the following more detailed form. See e.g.~\cite{BH,HH,MP,P}.

\begin{theorem}\label{converse} Let $(d_i)_{i \ge 0}$ be a sequence in $\N$ such that $d_0=1$ and $d_{i+1} \ \le \ d_i^{\vs{i}}$ for all $i \ge 1$.  Set $n=d_1$. There exists a lexideal $L$ in $S=K[X_1,\dots,X_n]$ such that for $R=S/L=\oplus_{i \ge 0} R_i$, we have $d_i=\dim R_i$ for all $i \ge 0$. 
\end{theorem}

This result is constructive, implying in turn that our results, namely Theorems~\ref{thm almost} and~\ref{thm main}, are constructive as well. A concrete illustration is given in the extended Example~\ref{concrete ex} below. One key point is the following intimate link between lexsegments and the numerical operation $a \mapsto a^{\vs{i}}$. 
\begin{lemma}\label{lem link} Let $C \subset \M_i$ be a lexsegment such that $|\M_i \setminus C|=a$. Then $|\M_{i+1} \setminus \M_1C | = a^{\vs{i}}$.
\end{lemma}

\subsection{An additive version of Macaulay's theorem}

Consider the abelian semigroup $G = \N^n$. For $1 \le i \le n$, denote by $e_i$ the $i$th canonical basis element of $G$, i.e. $e_i=(\delta_{ij})_{1 \le j \le n}$ where $\delta_{ij}$ is the Kronecker symbol. Let $B=\{e_1,\dots,e_n\} \subset G$. Note that for all $h \ge 1$, the $h$-fold iterated sumset $hB$ consists of all elements in $G$ whose coordinate sum is equal to $h$. Of course, $G$ is canonically isomorphic to the set $\M$ of monomials in $K[X_1,\dots,X_n]$, viewed as a multiplicative abelian semigroup. We order $G$ by transfering the graded lexicographic order $\le$ on $\M$ via the canonical isomorphism induced by $X_j \leftrightarrow e_j$ for all $j$. The following statement is equivalent to Macaulay's Theorem~\ref{thm macaulay}.

\begin{theorem}\label{thm additive macaulay} Let $G=\N^n$ and $B=\{e_1,\dots,e_n\} \subset G$. For all $h \ge 1$ and all subsets $A \subseteq hB$, we have
$$
|A+B| \ge |A^{\text{lex}}+B|,
$$
where $A^{\text{lex}} \subseteq hB$ denotes the unique lexsegment of cardinality $|A^{\text{lex}}|=|A|$. 
\end{theorem}
\begin{proof} See~\cite[Theorem 4.1]{MP} for an analogous statement in terms of monomial subspaces, shown there to be equivalent to Theorem~\ref{thm macaulay}.
\end{proof}

Macaulay's theorem is fundamental in commutative algebra and algebraic geometry, and since the 1970's in combinatorics too, thanks to the pioneering work on polytopes by McMullen~\cite{MC} and Stanley~\cite{S} among others. The additive version given by Theorem~\ref{thm additive macaulay} shows that Macaulay's theorem squarely belongs to additive combinatorics as well.

\section{An Almost Sharp Realization}\label{sec almost}

We show here that if $(d_i)_{i \ge 0}$ is a sequence of positive integers satisfying $d_0=1$ and
\begin{equation}\label{mac conditions}
1 \le d_{i+1} \le d_i^{\vs{i}} 
\end{equation}
for all $i \ge 1$, then there exists an abelian semigroup $G$ and a subset $A \subseteq G$ such that 
\begin{equation}\label{almost sharp}
d_h \le |hA| \le d_h+1
\end{equation}
for all $h \ge 0$. Our proof of this almost sharp realization relies on the sufficiency condition in Macaulay's theorem, and more specifically on Theorem~\ref{converse}. To proceed, we need a few relevant facts concerning monomial ideals.

\subsection{Monomial ideals}\label{sec monomial ideals}

Let $S=K[X_1,\dots,X_n]$ be the $n$-variable polynomial algebra over the field $K$, endowed with its standard grading $S=\oplus_{i \ge 0}S_i$ induced by $\deg(X_j)=1$ for all $j$. As earlier, we denote by $\M$ the set of monomials in $S$ and by $\M_i = \M \cap S_i$ the subset of monomials of degree $i$ for all $i \ge 0$.

A \emph{monomial ideal} in $S=K[X_1,\dots,X_n]$ is an ideal $J$ of $S$ generated by monomials. Of course, $J$ is a graded ideal, so that $J=\oplus_{i \ge 0} J_i$, where $J_i=J \cap S_i$. Macaulay proved that for every graded ideal $I$ of $S$, there exists a monomial ideal $J$ of $S$ such that $S/I$ and $S/J$ have the same Hilbert function. See Proposition~\ref{prop same hilbert function} below.

\begin{lemma}\label{basic lemma} Let $J \subset S$ be a monomial ideal. Let $f \in S$. Then $f \in J$ if and only if every monomial with a nonzero coefficient in $f$ belongs to $J$.
\end{lemma}

\begin{proof}
Easily follows from the fact that $J$ is spanned by monomials in $\M$ and that $\M$ is a $K$-basis of $S$.
\end{proof}

\begin{proposition}[Macaulay, \cite{M}]\label{monomial basis} Let $J \subset S$ be a monomial ideal. Let $\pi \colon S \to S/J$ be the quotient map. Then the family $\F=\{\pi(u) \mid u \in \M \setminus J\}$ is a $K$-basis of $S/J$. 
\end{proposition}

\begin{proof}
The family $\F$ spans $S/J$, since $\M$ spans $S$ and $\pi(\M \cap J)=\{0\}$. And $\F$ is free, for if $f=\sum_{u \in \M \setminus J} \lambda_u u$ and $\pi(f)=0$, then $f \in \ker(\pi)=J$. Lemma~\ref{basic lemma} then implies $\lambda_u=0$ for all $u \in \M \setminus J$, i.e. $f=0$. 
\end{proof}

Even though we have already encountered iterated product sets above, we formally recall the notation here.
\begin{notation}
Let $G$ be an abelian semigroup in multiplicative notation. For any subset $A \subseteq G$, we denote by $A^h=\underbrace{A \cdots A}_h$
its $h$-fold iterated product set.
\end{notation}

We need one more auxiliary result. 

\begin{proposition}\label{inter} Let $J$ be a monomial ideal in $S$. Let $R=S/J=\oplus_{i \ge 0} R_i$ and let
$
\pi \colon S \to R= S/J
$
be the quotient map. Let $x_j=\pi(X_j)$ for all $j$ and set $A=\{x_1,\dots,x_n\} \subset R$. Then for all $h \ge 1$, we have
\begin{equation}\label{formula inter}
|A^h| = 
\begin{cases}
\dim R_h & \textrm{if } J_h = \{0\}, \\
\dim R_h+1 & \textrm{if not},
\end{cases}
\end{equation}
where $J_i=S_i \cap J$ for all $i$. 
\end{proposition}
\begin{proof} We have $J = \oplus_{i \ge 0}J_i$, and $J_i$ has for vector subspace basis $\M_i \cap J$ for all $i \ge 0$. Since $A=\pi(\M_1)$, and since $\M_h=\M_1^h$ for all $h \ge 1$, we have 
\begin{equation}\label{A^h}
A^h=\pi(\M_h)
\end{equation} for all $h \ge 1$. Since $J = \ker(\pi)$, we have
\begin{equation}\label{pi M_h}
\pi(\M_h) = 
\begin{cases}
\pi(\M_h \setminus J_h) & \textrm{if }  \M_h \cap J_h = \emptyset, \\
\pi(\M_h \setminus J_h) \sqcup \{0\} & \textrm{if not.}
\end{cases}
\end{equation}
It follows from Proposition~\ref{monomial basis} that
\begin{equation}\label{basis M-J}
\dim R_h = |\M_h \setminus J_h|=|\pi(\M_h \setminus J_h)|.
\end{equation}
Combining~\eqref{A^h},~\eqref{pi M_h} and \eqref{basis M-J} yields the claimed formula~\eqref{formula inter}.
\end{proof}

\subsection{First construction}

Combining the above results with the sufficiency part of Macaulay's theorem, we obtain an almost sharp realization of $d_h$ as $|hA|$ for some subset $A$ of some abelian semigroup.

\begin{theorem}\label{thm almost} Let $(d_i)_{i \ge 0}$ be a sequence of nonnegative integers such that $d_0=1$ and $d_{i+1} \le d_i^{\vs{i}}$ for all $i \ge 1$. Then there exists an abelian semigroup $G$ and $A \subseteq G$ such that $d_h \le |hA| \le d_h+1$ for all $h \ge 0$.
\end{theorem}
\begin{proof} Set $n=d_1$. By Macaulay's theorem, there exists a standard graded algebra $R = \oplus_{i \ge 0} R_i$ such that $d_i=\dim R_i$ for all $i \ge 0$. By Theorem~\ref{converse}, one may take $R=S/L$, where $S=K[X_1,\dots,X_n]$ with its standard grading, and $L$ is a suitable lexideal in $S$. Let $\pi \colon S \to R$ be the quotient map. For the required abelian semigroup, in multiplicative notation, we may take $G=(R,\cdot)$ or, more economically, $G=\pi(\M)$. Set $x_j=\pi(X_j)$ for all $j$ and $A=\{x_1,\dots,x_n\} \subset G$. It then follows from Proposition~\ref{inter} that
$|A^h| \in \{d_h,d_h+1\}$ for all $h \ge 0$, as desired.
\end{proof}

Given a sequence $(d_i)_{i \ge 0}$ satisfying the conditions of Theorem~\ref{thm almost}, the following extended example shows how to explicitly construct a pair $G,A$ satisfying the conclusion of this theorem.

\begin{example}\label{concrete ex} Let $(d_0,d_1,d_2,d_3,d_4,d_5,\dots)=(1,5,13,25,42,63,\dots)$. 
Then $d_{i+1} \le d_i^{\vs{i}}$ for $1 \le i \le 4$. Indeed, we have
\begin{align*}
d_1 = 5 = \binom{5}{1} & \, \Longrightarrow\,  d_1^{\vs{1}} = \binom{6}{2}=15; \\
d_2 = 13 = \binom{5}{2}+\binom{3}{1} & \, \Longrightarrow\,  d_2^{\vs{2}}= \binom{6}{3}+ \binom{4}{2} = 26; \\
d_3 = 25 = \binom{6}{3}+\binom{3}{2}+\binom{2}{1} & \, \Longrightarrow\,  d_3^{\vs{3}}= \binom{7}{4}+ \binom{4}{3}+ \binom{3}{2} = 42; \\
d_4 = 42 = \binom{7}{4}+ \binom{4}{3}+ \binom{3}{2} & \, \Longrightarrow\, d_4^{\vs{4}} = \binom{8}{5}+ \binom{5}{4}+ \binom{4}{3}=65.
\end{align*}
Hence the differences $d_{i}^{\vs{i}}-d_{i+1}$ assume the following nonnegative values, as claimed: 
\begin{equation}\label{differences}
d_1^{\vs{1}}-d_2 = 2, \quad d_2^{\vs{2}}-d_3 = 1, \quad d_3^{\vs{3}}-d_4 = 0, \quad d_4^{\vs{4}}-d_5 = 2.
\end{equation}
Set $n=d_1=5$ and $S=K[X_1,\dots,X_5]$. We now use~\eqref{differences} to construct a lexideal $L \subset S$ such that the quotient $R=S/L=\oplus_{i \ge 0}R_i$ satisfies $\dim R_i = d_i$ for $0 \le i \le 5$. To do so, it suffices to exhibit a minimal system of monomial generators $\G$ satisfying the following requirements:
\begin{enumerate}
\item[(1)] $|\G \cap \M_{i+1}|=d_{i}^{\vs{i}}-d_{i+1}$ for all $1 \le i \le 4$,
\item[(2)] the resulting ideal $L=(\G)$ is a lexideal.
\end{enumerate}
The first condition arises from Lemma~\ref{lem link}. Using these constraints as a construction tool, we obtain the following solution: 
$$
\G=\{X_1^2,X_1X_2, X_1X_3^2, X_1X_3X_4^3,X_1X_3X_4^2X_5\}.
$$
As required, we do have $|\G \cap \M_{2}|=d_{1}^{\vs{1}}-d_{2}=2$, $|\G \cap \M_{3}|=1$, $|\G \cap \M_{4}|=0$, $|\G \cap \M_{5}|=2$, and $L \cap \M_i$ is a lexsegment for all $i \ge 2$. Let $\pi \colon S \to R=S/L$ be the quotient map. Again, the sought-for semigroup may be taken as $G=(R,\cdot)$, or more simply $G = \pi(\M)$. Set $x_j=\pi(X_j)$ for $1 \le j \le 5$, and $A = \{x_1,\dots,x_5\} \subset G$. Then
$$
|A|=d_1, \quad |A^h| = d_h+1
$$
for all $2 \le h \le 5$, as desired. For instance, for $h=2$ we have $x_1^2=x_1x_2=0$ in $G$, and
$$A^2 = \{0\} \sqcup \{x_1x_3, x_1x_4, x_1x_5, x_2^2, x_2x_3, x_2x_4, x_2x_5, x_3^2, x_3x_4, x_3x_5, x_4^2, x_4x_5, x_5^2 \},$$
so that $|A^2|=14=d_2+1$.
\end{example}

\section{Main Result}\label{sec full}
In order to show that the bounds given by Theorem~\ref{thm EM} are best possible,  we now aim for a sharp realization. That is, given any sequence $(d_i)_{i \ge 0}$ of positive integers satisfying $d_0=1$ and
\begin{equation*}
1 \le d_{i+1} \le d_i^{\vs{i}} 
\end{equation*}
for all $i \ge 1$, we shall construct an abelian semigroup $G$ and a subset $A \subseteq G$ such that 
\begin{equation*}
d_h = |hA|
\end{equation*}
for all $h \ge 0$. Note that the condition $d_h \ge 1$ for all $h$ is necessary here, since $|hA| \ge 1$ for any nonempty subset $A$ of a semigroup $(G,+)$. To that end, we shall deform the lexideal $L \subset S$, used above for our almost sharp realization, into a binomial ideal $\hat{L} \subset S$ with the same Hilbert function as $L$, i.e. such that $\dim \hat{L} \cap S_i=\dim L \cap S_i=$ for all $i$. The latter constraint can be achieved with a Gr\"{o}bner basis construction.

\subsection{Gr\"{o}bner bases}\label{sec grobner}
We recall here the few relevant facts on Gr\"obner bases needed for our constructions, and refer to~\cite{Fr, HH, P} for more details. Again, let $\M$ denote the set of monomials in $K[X_1, \dots, X_n]$. The notion of Gr\"{o}bner basis is relative to a given ordering of $\M$. Here we only consider the graded lexicographic ordering $\le$ on $\M$ relative to $X_1 > \dots > X_n$ as defined in Section~\ref{sec lexideals}. 

\smallskip
Denote $\M^* = \M \setminus \{1\}$. For any $u,v \in \M$, let $\gcd(u,v) \in \M$ denote their greatest common divisor. We further need the following notation and definitions.

\begin{notation}
For a nonzero polynomial $f \in K[X_1, \dots, X_n]$, we denote by $\inm(f) \in \M$ its \emph{leading monomial} with respect to the given ordering on $\M$, and by $\lc(f) \in K^*$ its \emph{leading coefficient}, i.e. the coefficient of $\inm(f)$ in $f$. The \emph{leading term} of $f$ is 
$$
\lt(f)=\lc(f)\inm(f).
$$
For a proper ideal $I \subsetneq K[X_1, \dots, X_n]$, we denote by $\inm(I)$ the monomial ideal generated by the set $\{\inm(f) \mid f \in I \setminus \{0\}\}$.
\end{notation}

The importance of the ideal $\inm(I)$ stems from the following property.

\begin{proposition}[Macaulay, \cite{M}]\label{prop same hilbert function} Let $I$ be a proper graded ideal in $S=K[X_1, \dots, X_n]$. Then the graded algebras
$$
S/I \text{ and } S/\inm(I)
$$
have the same Hilbert function.
\end{proposition}

\begin{definition} A finite set $\{g_1,\dots,g_s\} \subset K[X_1, \dots, X_n]\setminus K$ of nonconstant polynomials is a \emph{Gr\"{o}bner basis} if, for any nonzero element $f$ of the ideal $I=(g_1,\dots,g_s)$, we have $\inm(f) \in (\inm(g_1), \dots, \inm(g_s))$; equivalently, $\inm(f)$ is divisible by $\inm(g_i)$ for some $1 \le i \le s$. We then say that $\{g_1,\dots,g_s\}$ is a Gr\"obner basis of $I$. 
\end{definition}

Note that every proper ideal $I \subsetneq K[X_1, \dots, X_n]$ admits a Gr\"obner basis; this follows from the fact that $K[X_1, \dots, X_n]$ is noetherian, whence $\inm(I)$ is finitely generated. A key property of Gr\"{o}bner bases is the following direct consequence of Proposition~\ref{prop same hilbert function}. 

\begin{corollary}\label{cor same hilbert function} Let $\{g_1,\dots,g_s\} \subset K[X_1, \dots, X_n] \setminus K$ be a Gr\"{o}bner basis, with $g_j$ homogeneous for all $j$. Then the graded algebras
$$
S/(g_1,\dots,g_s) \textrm{ and } S/(\inm(g_1),\dots,\inm(g_s))
$$
have the same Hilbert function.
\end{corollary}
\begin{proof} Let $I=(g_1,\dots,g_s)$. Then $I$ is a graded ideal since the $g_j$ are homogeneous for all $j$. Moreover, $\inm(I)=(\inm(g_1),\dots,\inm(g_s))$ since $\{g_1,\dots,g_s\}$ is a Gr\"obner basis by hypothesis. We conclude with Proposition~\ref{prop same hilbert function}. 
\end{proof}

Buchberger developed an algorithm to construct Gr\"{o}bner bases for any proper ideal of $K[X_1,\dots,X_n]$, including a stopping criterion to recognize them. Here are the relevant details for the sequel.

\begin{definition}
Let $f,g,h \in K[X_1, \dots, X_n]$ with $f,h$ nonzero. We say that \emph{$f$ properly reduces to $g$ with respect to $h$} if $\inm(h)$ divides $\inm(f)$ in $\M$, and if $g$ is obtained by eliminating the leading term of $f$ with that of $h$, i.e.
$$
g=f-\frac{\lt(f)}{\lt(h)}h.
$$
We write $f \stackrel{h}{\longrightarrow} g$ when this occurs. In particular, if $f \stackrel{h}{\longrightarrow} g$, then either $g=0$ or else $\inm(g) < \inm(f)$. 
\end{definition}

\begin{definition}
More generally, let $H \subset K[X_1, \dots, X_n]$ be a set of nonconstant polynomials, and let $f,g \in K[X_1, \dots, X_n]$ with $f \neq 0$. We say that \emph{$f$ properly reduces to $g$ with respect to $H$}, and we write $f \stackrel{H}{\longrightarrow} g$, if there is a sequence of proper reductions from $f$ to $g$ of the form
$$
f =f_0 \stackrel{h_1}{\longrightarrow} f_1 \stackrel{h_2}{\longrightarrow} \dots \stackrel{h_\ell}{\longrightarrow} f_\ell = g
$$
with $h_1,\dots,h_\ell \in H$. 
\end{definition}

A key ingredient in Buchberger's algorithm is the notion of $\emph{$S$-polynomial}$.

\begin{definition}
Let $f,g \in K[X_1,\dots,X_n]\setminus K$. Let $v=\gcd(\inm(f),\inm(g)) \in \M$. The $\emph{$S$-polynomial}$ of $f,g$ is
$$
S(f,g) = \frac{\lt(g)}{v}f - \frac{\lt(f)}{v}g.
$$
\end{definition}

\begin{theorem}[Buchberger's criterion]\label{thm buchberger}
A set $H=\{f_1,\dots,f_r\}$ of polynomials in $K[X_1,\dots,X_n] \setminus K$ is a Gr\"{o}bner basis if and only if $S(f_i,f_j) \stackrel{H}{\longrightarrow} 0$ for all $1 \le i < j \le r$.
\end{theorem}

\subsection{A Gr\"{o}bner basis of binomials}
We construct here a Gr\"{o}bner basis made of certain homogeneous binomials, i.e. of differences $u-v$ of monomials $u,v$ of same degree. As above, $\M$ is the set of monomials in $K[X_1,\dots,X_n]$, endowed with the graded lexicographic order $\le$, and $\M^*=\M \setminus \{1\}$.

\begin{notation} For $u \in \M^*$, we denote by $\min(u)$ the smallest index $i \ge 1$ such that $X_i$ divides $u$, and by $\max(u)$ the largest index $j \ge 1$ such that $X_j$ divides $u$.
\end{notation}

For instance, for $u=X_2^4X_3X_5^3$, we have $\min(u)=2$ and $\max(u)=5$.

\begin{definition}\label{map phi}
Let $\varphi \colon \M^* \to \M^*$ be the map defined for all $u \in \M^*$ by 
$$
\varphi(u)=u{X_n}/{X_{\min(u)}}.
$$
\end{definition}
Note that if $\min(u) < n$, then $u > \varphi(u)$ and hence $\inm(u-\varphi(u))=u$. For instance, for $u=X_2^4X_3X_5^3$ again, taken here as an element of $K[X_1,\dots,X_5]$, i.e. with $n=5$, we have 
$$
\varphi(u)=X_2^3X_3X_5^4
$$
and, as stated, $X_2^4X_3X_5^3 > X_2^3X_3X_5^4$ in $\M_8$.

\begin{proposition}\label{binoms grobner} Let $u_1,\dots,u_r \in \M^*$ satisfy $\min (u_i) \le n-1$ for all $i$. Then the set of binomials
$$
H_r=\{u_i - \varphi(u_i) \mid 1 \le i \le r\}
$$
is a Gr\"{o}bner basis.
\end{proposition}
\begin{proof} The case $r=1$ is trivial. Let $r=2$, and let $u_1,u_2 \in \M^*$ satisfy $\min(u_1), \min(u_2) \le n-1$. With Theorem~\ref{thm buchberger} in mind, we will show that
\begin{equation}\label{s-poly}
S(u_1-\varphi(u_1), u_2-\varphi(u_2)) \stackrel{H_2}{\longrightarrow} 0.
\end{equation}
Without loss of generality, we may assume $u_1 > u_2$ and $\min(u_1)=1$. Let $i = \min(u_2)$. Thus $i \in \{1,\dots,n-1\}$ by hypothesis. Write $u_1=X_1v_1$ and $u_2=X_iv_2$ with $v_1,v_2 \in \M$ and $\min(v_1) \ge 1$, $\min(v_2) \ge i$. Then 
\begin{eqnarray*}
u_1-\varphi(u_1) & = & (X_1-X_n)v_1, \\ 
u_2-\varphi(u_2) & = & (X_i-X_n)v_2.
\end{eqnarray*}
Let now $v = \gcd(v_1,v_2) \in \M$.

\smallskip
$\bullet$ Assume first $i=1$. Then $X_1v=\gcd(u_1,u_2)$, and we have

\begin{align*}
S(u_1-\varphi(u_1),u_2-\varphi(u_2)) & =  S((X_1-X_n)v_1,(X_1-X_n)v_2) \\
& =  (X_1-X_n)v_1{v_2}/{v}-(X_1-X_n)v_2{v_1}/{v} \\
& = 0.
\end{align*}

$\bullet$ Assume now $i \ge 2$. Then

\begin{align*}
S(u_1-\varphi(u_1),u_2-\varphi(u_2)) & =  S((X_1-X_n)v_1,(X_i-X_n)v_2) \\
& =  (X_1-X_n)X_iv_1{v_2}/{v}-(X_i-X_n)X_1v_2{v_1}/{v} \\
& \stackrel{u_1-\varphi(u_1)}{\longrightarrow} X_1X_nv_2{v_1}/{v} - X_iX_nv_1{v_2}/{v} \\
& \stackrel{u_2-\varphi(u_2)}{\longrightarrow} X_n^2v_2{v_1}/{v} - X_n^2v_1{v_2}/{v}\\
& =  0.
\end{align*}

By Buchberger's criterion in Theorem~\ref{thm buchberger}, the set $H_2$ is a Gr\"{o}bner basis, as desired. For $r \ge 3$, the analog of formula~\eqref{s-poly} remains valid for any pair $u_i-\varphi(u_i),u_j-\varphi(u_j)$ with $1 \le i < j \le r$. Hence, by Buchberger's criterion again, the set $H_r$ is a Gr\"{o}bner basis, and the proof is complete.
\end{proof}

\subsection{Sharp realization}
We are now in a position to prove our main result in this paper.

\begin{theorem}\label{thm main} Let $(d_i)_{i \ge 0}$ be a sequence of positive integers such that $d_0=1$ and $1 \le d_{i+1} \le d_i^{\vs{i}}$ for all $i \ge 1$. Then there exists an abelian semigroup $G$ and $A \subseteq G$ such that $d_h = |hA|$ for all $h \ge 0$.
\end{theorem}
\begin{proof}
Set $n=d_1$ and $S=K[X_1,\dots,X_n] = \oplus_{i \ge 0} S_i$ with its standard grading. By Theorem~\ref{converse}, there exists a lexideal $L \subseteq S$ such that, for $R=S/L=\oplus_{i \ge 0}R_i$, we have
$$
d_i=\dim R_i
$$
for all $i \ge 0$. Denoting $L_i=L \cap S_i$, we have $L=\oplus_{i \ge 0}L_i$ and $R_i=S_i/L_i$ for all $i$. In particular, since $n=d_1=\dim R_1$, we have $L_1=\{0\}$.

\medskip
\noindent
\textbf{Claim 1.} For all $u \in \M \cap L$, we have
\begin{equation}\label{min u}
\min(u) \le n-1.
\end{equation}
For otherwise, let $u \in \M \cap L$ be such that $\min(u)=n$. Therefore $u=X_n^k$ for some $k \ge 1$. This implies $L \supset \M_k$ since $X_n^k=\min(\M_k)$ and $L \cap \M_k$ is a lexsegment. Hence $R_k = \{0\}$, contradicting $d_k \ge 1$. This proves the claim. 

\medskip
Let $\G \subset \M \cap L$ be the minimal system of monomial generators of $L$, so that $L=(\G)$. Of course $\G$ is finite and consists of all monomials in $L$ which are not the product of two monomials in $L$. Denote $\G_i=\G \cap S_i=\G \cap \M_i$ for all $i \ge 1$. We have $\G_1=\emptyset$ since $L_1=\{0\}$. 

\smallskip
Since $L$ is a lexideal, it is a \emph{stable} monomial ideal. As such, its minimal system of generators $\G$ may be characterized as follows~\cite{EK}: for all $u \in \M \cap L$, there is a unique monomial factorisation
$$
u = vw
$$
with $v,w \in \M$ such that
\begin{equation}\label{canonfac}
\begin{cases}
v \in \G, & \\
\max(v)  \le \min(w). &
\end{cases}
\end{equation}

Using the map $\varphi \colon \M^* \to \M^*$ of Definition~\ref{map phi}, denote
$$
\hat{\G}=\{u-\varphi(u) \mid u \in \G\}.
$$
It follows from~\eqref{min u} and the definition of $\varphi$ that $u > \varphi(u)$ for all $u \in \G$. Set $\hat{L} = (\hat{\G})$, the homogeneous binomial ideal of $S$ generated by $\hat{\G}$. Denote by 
\begin{align*}
\pi &\colon S \to R=S/L = \oplus_{i \ge 0} R_i, \\
\hat{\pi} &\colon S \to \hat{R}=S/\hat{L} =\oplus_{i \ge 0} \hat{R}_i
\end{align*}
the respective quotients and quotient maps of $S$. Applying Proposition~\ref{binoms grobner} to $\hat{\G}$, as allowed by~\eqref{min u}, it follows that $\hat{\G}$ is a Gr\"{o}bner basis of $\hat{L}$. Therefore, by Corollary~\ref{cor same hilbert function}, the Hilbert functions of $L$ and $\hat{L}$ are the same. That is, for all $i \ge 0$, we have
\begin{equation}\label{same dim}
\dim (\hat{R}_i) = \dim (R_i) =d_i.
\end{equation}

\medskip
\noindent
\textbf{Claim 2.} For all $u \in \M \cap L$, we have
\begin{equation}\label{congruence}
\varphi(u) \equiv u \bmod \hat{L}.
\end{equation}
Indeed, consider the unique monomial decomposition $u=vw$ with $v \in \G$ provided by~\eqref{canonfac}. Hence $\min(u)=\min(v)$, implying
$$
\varphi(u)=\varphi(vw)=\varphi(v)w
$$
by definition of $\varphi$. Therefore
$$
u-\varphi(u)=vw-\varphi(vw)=vw-\varphi(v)w=(v-\varphi(v))w.
$$
Since $v-\varphi(v) \in \hat{\G}$, this proves~\eqref{congruence}.

\medskip
\noindent
\textbf{Claim 3.} For all  $u \in \M \cap L$, there is a least exponent $\ell \ge 1$ such that 
\begin{equation}\label{least l}
\varphi^{\ell}(u) \in \M \setminus L,
\end{equation}
where $\varphi^{\ell}=\underbrace{\varphi \circ \dots \circ \varphi}_{\ell}$. Indeed, at each application of $\varphi$, the exponent of $X_n$ increases by $1$ while the degree remains constant. And $X_n^k \in \M \setminus L$ for all $k$ by~\eqref{min u}. This proves the claim.

\medskip
\noindent
\textbf{Claim 4.} We have $\hat{\pi}(\M)=\hat{\pi}(\M \setminus L).$ That is, for all $u \in \M$, there exists $w \in \M \setminus L$ such that
\begin{equation}\label{u equiv w}
u \equiv w \bmod \hat{L}.
\end{equation}
Indeed, if $u \in \M \setminus L$, take $w=u$. If $u \in \M \cap L$, let $w=\varphi^{\ell}(u) \in \M \setminus L$ with $\ell$ minimal as given by~\eqref{least l}. We have
$$
u-\varphi^{\ell}(u)=\sum_{i=0}^{\ell-1} (\varphi^{i}(u)-\varphi^{i+1}(u)).
$$
By minimality of $\ell$ with respect to~\eqref{least l}, we have $\varphi^i(u) \in \M \cap L$ for all $0 \le i \le \ell-1$. Hence, since $v \equiv \varphi(v) \bmod \hat{L}$ for all $v \in \M \cap L$ by~\eqref{congruence}, it follows that
$$
u \equiv \varphi^{\ell}(u) \bmod \hat{L}.
$$
This shows that~\eqref{u equiv w} holds with $w=\varphi^{\ell}(u) \in \M \setminus L$, as desired. This settles the claim.

\medskip
Finally, let $A=\hatpi{(\M_1)} \subset \hat{R}$. Then for all $h \ge 0$,  we have
\begin{equation}\label{A^h}
A^h =\hatpi(\M_h)=\hatpi(\M_h \setminus L_h)
\end{equation}
since $\M_h=\M_1^h$. Moreover, $\hatpi(\M \setminus L)$ is a $K$-basis of $\hat{R}$. This follows from the facts that $\hatpi(\M \setminus L)$ spans $\hat{R}$, that $\pi(\M \setminus L)$ is a $K$-basis of $R$ by Proposition~\ref{monomial basis}, and by~\eqref{same dim}. We conclude that
$$
|A^h|=\dim(\hat{R}_h)=\dim(R_h)=d_h
$$
for all $h \ge 0$, as desired. 
\end{proof}

\begin{example}\label{concrete ex contd} Revisiting Example~\ref{concrete ex}, let $(d_0,d_1,d_2,d_3,d_4,d_5,\dots)=(1,5,13,25,42,63,\dots)$. Let $S=K[X_1,\dots,X_5]$, and let $L \subset S$ be the lexideal with minimal monomial generating set 
$$
\G=\{X_1^2,X_1X_2, X_1X_3^2, X_1X_3X_4^3,X_1X_3X_4^2X_5\}.
$$
Using the map $\varphi$ from Definition~\ref{map phi}, let $\hat{\G}=\{u-\varphi(u) \mid u \in \G\}$. Then
$$
\hat{\G}=\{X_1^2-X_1X_5,X_1X_2-X_2X_5, X_1X_3^2-X_3^2X_5, X_1X_3X_4^3-X_3X_4^3X_5,X_1X_3X_4^2X_5-X_3X_4^2X_5^2\}.
$$
This is a Gr\"obner basis by Proposition~\ref{binoms grobner}. Let $\hat{L}=(\hat{\G})$. Denote by $\hatpi \colon S \to \hat{R}=S/\hat{L}$ the quotient map, set $x_j=\hatpi(X_j)$ for $1 \le j \le 5$, and $A = \{x_1,\dots,x_5\} \subset \hat{R}$. Then by Theorem~\ref{thm main}, we have
$$
|A^h| = d_h
$$
for all $0 \le h \le 5$, as desired. 
\end{example}

\section{Concluding Remarks}\label{sec conclusion}

Theorem~\ref{thm main} provides optimal upper bounds on the growth of iterated sumsets relative to all abelian semigroups. However, restricted to abelian \emph{groups} only, the analogous problem remains open. 

As recalled in Theorem~\ref{thm EM}, we have shown in~\cite{EM} that the arithmetic conditions $d_{i+1} \le d_i^{\vs{i}}$ for all $i$ are satisfied by all sequences $(d_i)_{i \ge 0}$ occuring as $d_h=|hA|$ for all $h$, with $A$ a nonempty finite subset of an abelian semigroup $G$. But the arithmetic conditions $d_{i+1} \le d_i^{\vs{i}}$ for all $i$ do not imply that this sequence is nondecreasing. Yet relative to groups, and in contrast with semigroups in general, monotonicity is a necessary condition, since 
\begin{equation}\label{monoton}
|hA| \le |(h+1)A|
\end{equation}
for all finite subsets $A$ of groups and for all $h$. On the other hand, this monotonicity condition is far from being sufficient. For instance, consider the eventually constant sequence
\begin{eqnarray*}\label{13344}
(d_0, d_1, d_2, d_3, \dots) = (1,3,3,4,4,4, \dots)
\end{eqnarray*}
with $d_i = 4$ for all $i \ge 3$. The conditions $d_{i+1} \le d_i^{\vs{i}}$ for all $i \ge 1$ are satisfied here. Yet this sequence cannot be of the form $(|hA|)_{h \ge 0}$ for a subset $A$ of a group $G$. For if $|A|=|2A|=3$, then $A$ is a translate of a subgroup of order $3$, whence $|hA|=3$ for all $h \ge 1$. This follows from the following well known lemma, whose short proof we recall for convenience.

\begin{lemma} Let $A$ be a nonempty finite subset of a group $(G,+)$. Then $|hA|=|(h+1)A|$ for some $h \ge 0$ if and only if $A$ is a translate of a subgroup of $G$ of cardinality $|hA|$. In particular, $|hA|=|h'A|$ for all $h' \ge h$.
\end{lemma}
\begin{proof} Without loss a generality, we may assume that $A$ contains $0$. Let $B=hA$, where $h \ge 0$ satisfies $|hA|=|(h+1)A|$. Then $|hA| = |h'A|$ for all $h' \ge h$, and in fact $hA=h'A$ for all $h' \ge h$ since $hA \subseteq h'A$ as $A$ contains $0$. It follows that $B$ is a finite subset of $G$ satisfying $0 \in B$ and $2B=B$. Hence, for all $b \in B$, there exists $c \in B$ such that $b+c=0$. It follows that $B$ is both stable under addition and taking opposites. Hence it is a subgroup of $G$.
\end{proof}

Besides the problem of finding an analog of Theorem~\ref{thm main} restricted to abelian groups, the following general problem is completely open.

\medskip
\noindent
\textbf{Problem.} Characterize the nondecreasing sequences $(d_i)_{i \ge 0}$ of positive integers arising as iterated sumsets cardinalities in abelian groups, i.e. such that there exists an abelian group $G$ and a nonempty finite subset $A \subseteq G$ such that $d_h=|hA|$ for all $h \ge 0$.

\medskip
\noindent
\textbf{Acknowledgement.} This research was supported in part by the International Centre for Theoretical Sciences (ICTS) in Bangalore during a visit for the program - Workshop on Additive Combinatorics (Code: ICTS/wac2020/02).

\nocite{*}

\bibliographystyle{plain}

\bibliography{optimal}

\medskip
\noindent
\textbf{Authors' addresses:}

\begin{itemize}
\item[$\triangleright$] Shalom Eliahou, \\
Univ. Littoral C\^ote d'Opale, UR 2597 - LMPA - Laboratoire de Math\'ematiques Pures et Appliqu\'ees Joseph Liouville, F-62100 Calais, France and CNRS, FR2037, France. \\
\texttt{eliahou@univ-littoral.fr}
\item[$\triangleright$] Eshita Mazumdar, \\
School of Arts and Sciences, Ahmedabad University, Central Campus, Navrangpura, Ahmedabad 380009, India.\\
\texttt{eshita.mazumdar@ahduni.edu.in}
\end{itemize}

\end{document}